\documentclass[12pt]{article}
\usepackage{amsfonts}
\usepackage{theorem}

\usepackage[OT2,OT1]{fontenc}
\newcommand\cyr{%
\renewcommand\rmdefault{wncyr}%
\renewcommand\sfdefault{wncyss}%
\renewcommand\encodingdefault{OT2}%
\normalfont
\selectfont}
\DeclareTextFontCommand{\textcyr}{\cyr}

\newtheorem{theorem}{Theorem}

\newtheorem{lemma}{Lemma}

\usepackage[russian, english]{babel}
\usepackage{centernot}

\usepackage[unicode,colorlinks,plainpages=false]{hyperref}

\newcommand{\openbox}{$\begin{array}{c}
\hspace*{-0.55em}\sqcap \hspace*{-0.60em}\\[-0.4em] \hline
\multicolumn{1}{c}{\hspace*{-0.60em}}\\[-0.8em]
\end{array}$}

\begin{document}

\centerline{\bf On Congruence Permutable G-sets\footnote{{\bf Keywords:} G-sets; congruence permutable algebras; semigroups.\\
{\bf Mathematics Subject Classification:} 20E15; 20M05. This work was supported by the National Research, Development and Innovation Office – NKFIH, 115288\\ e-mail: nagyat@math.bme.hu}}

\medskip
\centerline{Attila Nagy}

\medskip
\centerline{Department of Algebra}
\centerline{Budapest University of Technology and Economics}

\bigskip

\bigskip

{\bf Abstract} An algebraic structure is said to be congruence permutable if its arbitrary
congruences $\alpha$ and $\beta$ satisfy the equation $\alpha \circ \beta =\beta \circ \alpha$, where $\circ$ denotes
the usual composition of binary relations. For an arbitrary $G$-set $X$ with $G\cap X=\emptyset$, we define a semigroup $(G,X,0)$ with a zero $0$ ($0\notin G\cup X$), and give necessary and sufficient conditions for the congruence permutability of the $G$-set $X$ by the help of the semigroup $(G,X,0)$.

\bigskip

\section{Introduction}

Let $G$ be a group with the identity element $e$. By a $G$-set we shall mean a right $G$-set, that is, a non-empty set $X$ together with a mapping \[X\times G \mapsto X;\quad (x, g)\mapsto x^g\in X\] satisfying the equations
$x^e=x$ and $(x^g)^h=x^{(gh)}$
%\[\varphi (x, e)=x\] and \[\varphi (\varphi (x, g_1), g_2)=\varphi (x, g_1g_2)\]
for every $x\in X$ and every $g, h\in G$.

A $G$-set $X$ is said to be transitive if for every $x, y\in X$ there is a $g\in G$ such that $x^g=y$. A transitive $G$-subset of a $G$-set $X$ is called an orbit of $X$. Clearly, any $G$-set is a disjoint union of its orbits.

Every $G$-set $X$ can be considered as a unary algebra $(X; G)$ with the set $G$ of operations where the operation $g\in G$ is defined by the role $g(x)=x^g$ for every $x\in X$.

By a congruence of a $G$-set $X$ we mean an equivalence relation $\sigma$ of $X$ which satisfies the following condition: for every $a, b\in X$, the assumption $(a, b)\in \sigma$ implies $(a^g, b^g)\in \sigma$ for every $g\in G$ (that is, $\sigma$ is a congruence of the unary algebra $(X,G)$).

\medskip

The next lemma is about the congruence lattice of a transitive $G$-set $X$ (see \cite[Lemma 3]{Pálfy} and \cite[Lemma 4.20]{McKenzie}).

\begin{lemma}\label{lem1}
Let $X$ be a $G$-set such that the group $G$ acts on $X$ transitively. Then the congruence lattice $Con(X)$ of the $G$-set $X$ is isomorphic to the interval $[Stab_G(x), G]$ of the subgroup lattice of $G$, where $x$ is an arbitrary element of $X$ and $Stab_G(x)=\{ g\in G:\ x^g=x\}$.
The corresponding isomorphisms are
$\phi :\alpha \mapsto H_{\alpha}=\{ g\in G:\ (x^g, x)\in \alpha\}$
and
$\psi : H\ \mapsto \alpha _H=\{ (x^g,x^h)\in A\times A:\ Hg=Hh\}$
($\alpha \in {\it Con(X)}, \ H\in [Stab_G(x),G])$) which are inverses of each other.
\hfill\openbox
\end{lemma}

\medskip

An algebraic structure is said to be congruence permutable if $\alpha \circ \beta =\beta \circ \alpha$ is satisfied for its arbitrary congruences $\alpha$ and $\beta$.
\medskip

By \cite[Lemma 1]{Pálfy2}, $\alpha \circ \beta =\beta \circ \alpha$ is satisfied for congruences $\alpha$ and $\beta$ of a transitive $G$-set $X$ if and only if $H_{\alpha}H_{\beta}=H_{\beta}H_{\alpha}$ is satisfied. Thus the following lemma is a characterization of the congruence permutable transitive $G$-sets.

\begin{lemma}\label{lem2} A transitive $G$-set $X$ is congruence permutable if and only if $HK=KH$ is satisfied for every subgroups $H$ and $K$ of $G$ belonging to the interval $[Stab_G(x), G]$, where $x$ is an arbitrary element of $X$.\hfill\openbox
\end{lemma}

\medskip

Arbitrary congruence permutable $G$-sets are characterized in \cite{Vernikov}.
A $G$-set $X$ is called segregated if every congruence $\alpha$ of the $G$-set $X$ satisfies the following condition: if $A$ and $B$ are different orbits of $X$ such that $(a_0, b_0)\in \alpha$ for some $a_0\in A$ and $b_0\in B$ then $(a, b)\in \alpha$ for all $a, b\in A\cup B$. By \cite[Theorem 3.4]{Vernikov} the following lemma is true.

\begin{lemma}\label{lmarbitrary} A $G$-set $X$ is congruence permutable if and only if $X$ is a segregated $G$-set such that $X$ has at most two orbits and every orbit of $X$ is congruence permutable.
\end{lemma}

\medskip

In this paper we give a semigroup theoretical characterization of congruence permutable $G$-sets. For an arbitrary $G$-set $X$ with $G\cap X=\emptyset$, we define a semigroup $(G,X,0)$ with a zero $0$ ($0\notin G\cup X$), and give necessary and sufficient conditions for the congruence permutability of the $G$-set $X$ by the help of the semigroup $(G,X,0)$.

\medskip

For semigroup theoretical terminologies used in our investigation, we refer to the paper \cite{Hamilton} and the book \cite{Nagy}.

\section{Results}

It is clear that every $G$-set is isomorphic to a $G$-set $X$ with $G\cap X=\emptyset$. In the next we suppose that the considered $G$-sets $X$ satisfy this condition.

\medskip

{\bf Construction}
Let $X$ be a right $G$-set (with condition $G\cap X =\emptyset$). Let $0$ be a symbol with $0\notin G\cup X$. On the set $S=G\cup X\cup \{ 0\}$, define an operation $*$ as follows. For arbitrary $g, h\in G$, let $g * h =gh$, where $gh$ is the original product of $g$ and $h$ in $G$. For arbitrary $x\in X$ and arbitrary $g\in G$, let $x * g=x^g$. Let $0 * g=0$ for every $g\in G$. If $a\in X\cup \{ 0\}$ then, for arbitrary $s\in S$, let $s * a=0$. It is easy to check that $S$ is a semigroup in which $0$ is the zero element, $G$ is a subgroup of $S$, and $X\cup \{ 0\}$ is a zero subsemigroup of $S$ (that is, $a * b=0$ for all $a, b\in X\cup \{ 0\}$ ). The semigroup $S$ will be denoted by $(G,X,0)$.\hfill\openbox

\medskip

In this paper we give a necessary and sufficient condition for the congruence permutability of a $G$-set $X$ by the help of the semigroup $(G,X,0)$.

\medskip

The next example shows that the congruence permutability of a $G$-set $X$ and the congruence permutability of the semigroup $(G,X,0)$ are not equivalent conditions, in general.

\bigskip

\noindent
{\bf Example} Let $X=\{ a, b\}$ be a two-element set and $G$ be an arbitrary group. Assume $a^g=a$ and $b^g=b$ for every $g\in G$. Then the orbits of the $G$-set $X$ are $\{ a\}$ and $\{ b\}$. It is clear that $X$ is a congruence permutable $G$-set.
Let $\alpha$ and $\beta$ be equivalence relations on the semigroup $S=(G, X, 0)$ whose classes are
$\alpha: \{ a; 0\}, \{ b\}, G$ and $\beta : \{b; 0\}, \{ a\}, G$.
It is easy to see that $\alpha$ and $\beta$ are congruences on the semigroup $(G,X,0)$. Since $(a, 0)\in \alpha$ and $(0, b)\in \beta$, then we have
$(a, b)\in \alpha \circ \beta$. If the semigroup $(G,X,0)$ was congruence permutable then we would have $(a, b)\in \beta \circ \alpha$ from which we would get
$(a, t)\in \beta$ and $(t, b)\in \alpha$ for some $t\in (G,X,0)$. Since $[a]_{\beta }=\{ a\}$ and $[b]_{\alpha} =\{ b\}$, we would get $a=b$ which is a contradiction. Consequently the semigroup $(G,X,0)$ is not congruence permutable.

\bigskip
The next theorem characterizes the transitive congruence permutable $G$-sets by the help of the semigroup $(G,X,0)$.

\begin{theorem}\label{thm1} A $G$-set $X$ is transitive and congruence permutable if and only if the semigroup $S=(G,X,0)$ is congruence permutable.
\end{theorem}

\noindent
{\bf Proof}.
Assume that $X$ is a transitive congruence permutable $G$-set. Let $N$ denote the set $X\cup \{ 0\}$. First we show that, for an arbitrary non-universal congruence $\alpha$ on the semigroup $S=(G,X,0)$, we have $[g]_{\alpha}\subseteq G$ for every $g\in G$, and $[0]_{\alpha}=\{ 0\}$ or $[0]_{\alpha}=N$. Let $\alpha$ be a non-universal congruence on the semigroup $S=(G,X,0)$.
Assume $(a,g)\in \alpha$ for some $a\in N, g\in G$. Then $(e * a,g)\in \alpha$, where $e$ is the
identity element of $G$. As $e * a=0$, we get $g\in [0]_{\alpha}$ from which it follows that $G\subseteq [0]_{\alpha}$. Let $a\in X$ be an arbitrary element. Then $X=a * G\subseteq [0]_{\alpha}$
and so $[0]_{\alpha}=S$. This contradicts the assumption that $\alpha$ is a non-universal congruence on $S$. Consequently
$[a]_{\alpha}\subseteq N$ and $[g]_{\alpha}\subseteq G$ for every $a\in N$ and every $g\in G$. Consider the case when $[0]_{\alpha}\neq \{ 0\}$. Then there is an element $a\in X$ such that $a\in [0]_{\alpha}$ and so
$X=a * G\subseteq [0]_{\alpha}$. Hence $[0]_{\alpha}=N$.

Let $\alpha$ and $\beta$ be arbitrary congruences on the semigroup $S=(G,X,0)$. We show that
$\alpha \circ \beta =\beta \circ \alpha$. We can suppose that $\alpha$ and $\beta$
are not the universal relations of $S$.
Let $b, c\in S$ be arbitrary elements. Assume $(b,c)\in \alpha \circ \beta$. Then there
is an element $x\in S$ such that $(b,x)\in \alpha$ and $(x,c)\in \beta$.
We have two cases.

Case 1: $x\in G$. In this case $b,c\in G.$ As $G$ is congruence permutable,
there is an element $y\in G$ with $(b,y)\in \beta$ and $(y,c)\in \alpha$.
Hence $(b,c)\in \beta \circ \alpha$.

Case 2: $x\in N=X\cup \{ 0\}$. In this case $b,c\in N$. We have two subcases.
If $[0]_{\beta}=N$ or $[0]_{\alpha}=N$, then $(b,c)\in \beta \cup \alpha \subseteq \beta \circ \alpha$.
Consider the case $[0]_{\beta}=[0]_{\alpha}=\{ 0\}$. In this case $X$ is saturated
by both $\alpha$ and $\beta$.
If $x=0$, then $b=c=0$ and so
$(b,c)\in \beta \circ \alpha$.
If $x\in X$, then $b,c\in X.$ Let
$\alpha ^+$ and $\beta ^+$ denote the restriction of $\alpha$ and $\beta$ to $X$. Then $\alpha ^+$ and $\beta ^+$ are congruences on the $G$-set $X$. Moreover $(b, c)\in \alpha ^+\circ \beta ^+$.
Since $X$ is a congruence permutable $G$-set, we get $(b, c)\in \beta ^+\circ \alpha ^+$. Then there is an element $y\in X$ such that
$(b, y)\in \beta ^+$ and $(y, c)\in \alpha ^+$ from which we get $(b, y)\in \beta$ and $(y, c)\in \alpha$, that is, $(b, c)\in \beta \circ \alpha$.

Thus we have $(b,c)\in \beta \circ \alpha$ in both cases.
Consequently, $\alpha \circ \beta \subseteq \beta \circ \alpha$, and by symmetry
$\alpha \circ \beta =\beta \circ \alpha$. Thus $S=(G,X,0)$ is a congruence permutable semigroup.

To prove the converse assertion, assume that the semigroup $S=(G,X,0)$ is congruence permutable. Let $\alpha, \beta $ be arbitrary congruences of the $G$-set $X$.
Let $\alpha '$ be the equivalence relation on the semigroup $S=(G,X,0)$ defined by $\alpha '=\alpha \cup \iota _S$, where $\iota _S$ denotes the identity relation on $S$. We show that $\alpha '$ is a congruence relation on $S$. Assume $(a, b)\in \alpha '$ for some $a, b\in S$. We can suppose that $a\neq b$. Then $a, b\in X$ and $(a, b)\in \alpha$. Let $s\in S$ be an arbitrary elements. Since $s * a = 0 = s * b$, then $(s * a, s * b)\in \alpha '$, and so $\alpha '$ is a left congruence on the semigroup $S$. If $s\in G$, then $a * s= a^s$ and $\quad b * s=b^s$ and so $(a * s, b * s)\in \alpha \subseteq \alpha '$. If $s\in X\cup \{ 0\}$, then
$a * s=0=b * s$ and $(a * s, b * s)\in \alpha '$. Hence $\alpha '$ is a right congruence on $S$. Consequently $\alpha '$ is a congruence on $S$.
Similarly, $\beta '$ defined by $\beta '=\beta \cup \iota _S$ is a congruence on the semigroup $S=(G,X,0)$.
We show that $\alpha \circ \beta =\beta \circ \alpha$.
Let $a, b\in X$ be arbitrary elements. Assume $(a, b)\in \alpha \circ \beta$. Then there is an element $x\in X$ such that
$(a, x)\in \alpha$ and $(x, b)\in \beta$.
As $\alpha \subseteq \alpha '$ and $\beta \subseteq \beta '$, we have
$(a, b)\in \alpha '\circ \beta '$. Since $S=(G,X,0)$ is a congruence permutable semigroup, then $(a, b)\in \beta '\circ \alpha '$ and so there is an element $t\in S=(G,X,0)$ such that $(a, t)\in \beta '$ and $(t, b)\in \alpha '$. As $X$ is saturated by $\alpha '$ and $\beta '$, we have $t\in X$ and so
$(a, t)\in \beta$ and $(t, b)\in \alpha$. Hence $(a,b)\in \beta \circ \alpha$. Consequently
$\alpha \circ \beta \subseteq \beta \circ \alpha$, and by symmetry $\alpha \circ \beta =\beta \circ \alpha$.
Hence $X$ is a congruence permutable $G$-set.

Assume that $X$ has at least two orbits. Let $A$ and $B$ be different orbits of $X$. It is clear that $A\cup \{ 0\}$ and $B\cup \{ 0\}$ are ideals of the semigroup $(G,X,0)$. By \cite[Theorem 4]{Hamilton}, the ideals of a congruence permutable semigroup form a chain with respect to inclusion. Then $A\subseteq B$ or $B\subseteq A$ which contradicts $A\cap B=\emptyset$. Consequently $X$ has one orbit. Thus $X$ is a transitive congruence permutable $G$-set.
\hfill\openbox

\bigskip

Let $X$ be a $G$-set.
We say that the semigroup $(G,X,0)$ is segregated if the $G$-set $X$ is segregated.

\begin{lemma} Let $X$ be a $G$-set. Then the semigroup $(G,X,0)$ is segregated if and only if every congruence $\alpha$ on $(G,X,0)$ satisfies the following condition: if $A$ and $B$ are different orbits of $X$ such that $(a_0, b_0)\in \alpha$ for some $a_0\in A$ and $b_0\in B$ then $(a, b)\in \alpha$ for all $a, b\in A\cup B$.
\end{lemma}

\noindent
{\bf Proof}. It is clear that if $\alpha$ is a congruence on the semigroup $(G,X,0)$, then the restriction of $\alpha$ to $X$ is a congruence of the $G$-set $X$. Moreover, if $\alpha$ is a congruence of the $G$-set $X$, then $\alpha '=\alpha \cup \iota _S$ is a congruence on the semigroup $S=(G,X,0)$, where $\iota _S$ denotes the identity relation on $S=(G,X,0)$. Thus the assertion of the lemma is obvious.\hfill\openbox

\medskip

Let $A$ be an orbit of a $G$-set $X$. The subsemigroup $(G,A,0)$ is called an orbit subsemigroup of the semigroup $(G,X,0)$. The next theorem characterizes arbitrary $G$-sets by the help of the semigroup $(G,X,0)$ and the orbit subsemigroups of $(G,X,0)$.

\begin{theorem}\label{thm6} A $G$-set $X$ is congruence permutable if and only if the semigroup $(G,X,0)$ is segregated such that it has at most two orbit subsemigroups, and every orbit subsemigroup of $(G,X,0)$ is congruence permutable.
\end{theorem}

\noindent
{\bf Proof}. Let a $G$-set $X$ be congruence permutable. By Lemma~\ref{lmarbitrary}, $X$ is a segregated $G$-set such that $X$ has at most two orbits and every orbit of $X$ is a congruence permutable transitive $G$-set. Then the semigroup $(G,X,0)$ is segregated by definition, and it contains at most two orbit subsemigroups. By Theorem~\ref{thm1}, every orbit subsemigroup of $(G,X,0)$ is congruence permutable.

Conversely, assume that the semigroup $(G,X,0)$ is segregated such that it has at most two orbit subsemigroups, and every orbit subsemigroup of $(G,X,0)$ is congruence permutable. Then the $G$-set $X$ is segregated by definition, and it has at most two orbits. Every orbit of $X$ is a congruence permutable $G$-set by Theorem~\ref{thm1}. Consequently $X$ is a congruence permutable $G$-set by Lemma~\ref{lmarbitrary}.\hfill\openbox

\end{document}